\documentclass{amsart}

\usepackage[english]{babel}
\usepackage{leftidx}
\usepackage[numbers]{natbib}
\usepackage{bigints}
\usepackage{empheq}
\usepackage{array}
\usepackage{graphicx}
\usepackage{lipsum}
\usepackage{amsthm}
\usepackage{needspace}
\usepackage[makeroom]{cancel}
\newenvironment{myproof}[2] {\paragraph{\textbf{Proof of {#1} {#2} :}}}{\hfill$\square$}
\newtheorem{theorem}{Theorem}[section]
\newtheorem{proposition}[theorem]{Proposition}
\newtheorem{lema}[theorem]{Lemma}
\newtheorem{question-non}[]{}

\newtheorem{cor}[theorem]{Corollary}
\newtheorem{definition}[theorem]{Definition}
\newtheorem{observation}[theorem]{Remark}
\newtheorem{example}[theorem]{Example}

\title{PROPERTIES OF COMPLETE NONCOMPACT WARPED PRODUCT GRADIENT YAMABE SOLITONS}

\author{Tokura, W. $^{1}$}
\address{$^{1}$ Universidade Federal de Goi\'as, IME, 131, 74001-970, Goi\^ania, GO, Brazil.}
\email{williamisaotokura@hotmail.com $^{1}$}

\author{Adriano, L. $^{2}$}
\address{$^{2}$ Universidade Federal de Goi\'as, IME, 131, 74001-970, Goi\^ania, GO, Brazil.}
\email{levi@ufg.br $^{2}$}

\author{Pina, R. $^{3}$}
\address{$^{3}$ Universidade Federal de Goi\'as, IME, 131, 74001-970, Goi\^ania, GO, Brazil.}
\email{romildo@ufg.br $^{3}$}

\author{Barboza, M. $^{4}$}
\address{$^{4}$ Insituto Federal Goiano, 75790-000, Rodovia Geraldo Silva Nascimento Km 2,5, Uruta\'i, GO, Brazil.}
\email{marcelo.barboza@ifgoiano.edu.br $^{4}$}

\thanks{$^{1,4}$ Supported by CAPES}



\begin{document}

\begin{abstract}
In this paper, we look for properties of gradient Yamabe solitons on top of warped product manifolds. Utilizing the maximum principle, we find lower bound estimates for both the potential function of the soliton and the scalar curvature of the warped product. By slightly modifying Li-Yau's technique so that we can handle drifting Laplacians, we were able to find three different gradient estimates for the warping function, one for each sign of the scalar curvature of the fiber manifold. As an application, we exhibit a nonexistence theorem for gradient Yamabe solitons possessing certain metric properties on the base of the warped product.
\end{abstract}

\maketitle

\section{Introduction and main results}
\label{intro}

An ordered quadruple $(M^{n},g,X,\rho)$ built up from a Riemannian manifold $(M^{n},g)$, a smooth vector field $X$ on $M^{n}$ and some constant $\rho\in\mathbb{R}$ is said to be a \textit{Yamabe soliton} provided that
\begin{equation}\label{def1}
(scal_{g}-\rho)g=\frac{1}{2}\mathfrak{L}_{X}g,
\end{equation}
where $\mathfrak{L}_{X}g$ is the Lie derivative of $g$ in the direction of $X$ and $scal_{g}$ is the scalar curvature of $(M^{n},g)$. The soliton is classified into three types according to the sign of $\rho$: expanding if $\rho<0$, steady if $\rho=0$ and shrinking if $\rho>0$. It may happen that $X=\nabla_{g}h$ is the gradient field of a smooth real function $h$ on $M^{n}$, called \textit{potential}, in which case the soliton is referred to as a \textit{gradient Yamabe soliton}. Equation \eqref{def1} then becomes	\begin{equation}\label{def2}(scal_{g}-\rho)g=Hess_{g}h,
\end{equation}
where $Hess_{g}h$ is the Hessian of $h$. The gradient soliton is called \textit{trivial} if $h$ is a constant function.

Objects responding by the names of \textit{almost Yamabe soliton} and \textit{almost gradient Yamabe soliton} are obtained from the above equations, both \eqref{def1} and \eqref{def2}, respectively, if $\rho\in\mathbb{R}$ is replaced with a smooth function $\rho:M^{n}\rightarrow\mathbb{R}$.

\begin{definition}\label{def0}
	The warped product of Riemannian manifolds $(B^{n},g_{B})$ and $(F^{d},g_{F})$, whose warpage is measured by a smooth function $f>0$ on $B^{n}$, corresponds to the product space $\hat{M}^{n+d}=B^{n}\times F$ alongside the following Riemannian metric tensor
	\begin{equation}\label{warped}
	\hat{g}=\pi_{B}^{\ast}g_{B}+f(\pi_{B})^{2}\pi_{F}^{\ast}g_{F},
	\end{equation}
	where $\pi_{B}:M^{n+d}\rightarrow B^{n}$ and $\pi_{F}:M^{n+d}\rightarrow F^{d}$ are the projections of $M^{n+d}$ on its first and second factors, respectively.
\end{definition}

Warped product manifolds have already proven themselves to be a rich source of examples in a wide range of distinct geometrical objects, of which solitons are an example (cf. \cite{al2012warped, barboza2018invariant, de2017gradient, feitosa2017construction, kim2013warped, tokura2019warped}). A warped product manifold like the one in Definition \ref{def0} is suggestively written in the form $B^{n}\times_{f}F^{d}$ and its components $B$, $F$ and $f$ are called base, fiber and warping function, in this order.

In paper \cite{tokura2019warped} the present authors have studies warped product gradient Yamabe soliton with assumption
	\vspace{0,1cm}
	\begin{equation}\label{3}
	scal_{g_F}=\lambda_{F}=\text{constant}, \hspace{0,5cm} \widetilde{h}=h\circ\pi,\hspace{0,5cm}h\in C^{\infty}(B),
	\end{equation}
	where $\widetilde{h}$ is the potential function of $B\times_{f}F$. As pointed out in \cite{tokura2019warped}, the above conditions are natural and bring to light the fact  that the topology of the base space imposes constraints on the analysis of the warped product. Inspired by this work, in this paper we will continue investigating the warped product gradient Yamabe solitons satisfying \eqref{3}.

In recent years, much efforts have been devoted to understanding the geometry of gradient Yamabe solitons. Under a integral assumption on the Ricci curvature, Wu \cite{wu2011class} provides an estimate for the scalar curvature as a function of $\rho$. As noticed by Wu this result excludes the analysis of Einstein solitons with negative constant curvature. Observing such gap, Chu \cite{chu2013scalar} improve this result by considering a lower bound of Bakry-\'Emery Ricci tensor.

\begin{theorem}\label{chu}\textup{(\cite{chu2013scalar})} Let $(M^{n},g,\nabla h,\rho)$ be an $n$-dimensional complete noncompact gradient Yamabe soliton with 
	$$^{M}Ric+Hess h_{1}\geq K,\quad \text{where}\quad h_{1}=-\frac{h}{2(n-1)},$$
for some constant $K\in\mathbb{R}$ and consider $R_{\ast}=\inf_{M}scal_{g} $.
\begin{itemize}
	\item If $\rho>0$, then $0\leq R_{\ast}\leq\rho$.
	\item If $\rho=0$, then $R_{\ast}=0$.
	\item If $\rho<0$, then $\rho\leq R_{\ast}\leq 0$.
\end{itemize}
\end{theorem}
From section \ref{examples} below, we see that $\hat{M}=\mathbb{H}^{3}\times_{f}(\mathbb{S}^{3}\times\mathbb{H}^{3})$ (Example \ref{example1}) is a complete noncompact expanding gradient Yamabe soliton with constant scalar curvature $-6$. By a straightforward computation we obtain the following expression for Ricci tensor for fields on lift $\mathcal{L}(\{p\}\times \mathbb{H}^{3})$,
\begin{equation*}
^{\hat{M}}Ric(V,W)=-\left(2x_{3}^{-\tfrac{8}{7}}+\tfrac{40}{49}\right)\hat{g}(V,W),\qquad p\in\mathbb{S}^{3},
\end{equation*}
which is unbounded from below. Then Theorem \ref{chu} exclude this solitons.

Our first theorem improves this result by considering a lower bound for the Bakry-\'Emery Ricci tensor of the base.

\begin{theorem}\label{th1}
	Let $(\hat{M}=B^{n}\times_{f}F^{d},\hat{g},\nabla \tilde{h},\rho)$ be a complete gradient Yamabe soliton satisfying 
	$$^{B}Ric+Hess w\geq K,\quad\text{where}\quad w=-d\log f-\frac{h}{2(n+d-1)},$$
	for some constant $K\in \mathbb{R}$ and consider $R_{\ast}=\inf_{\hat{M}} scal_{\hat{g}}$, $\psi_{*}=\inf_{B}\psi$, where $\psi=scal_{g_B}-\lambda$ and
	
	\begin{equation*}
	\lambda=-\frac{\lambda_{F}}{f^2}+\frac{2d}{f}\Delta f+\frac{d(d-1)}{f^{2}}g_{B}(\nabla f,\nabla f)+\rho,
	\end{equation*}

	\begin{enumerate}
		\item If $\rho>0$, then $0\leq R_{*}\leq \rho$. Furthermore, if $\psi(x_{0})=\psi_{*}=-\rho$ for some $x_{0}\in B$, then $\psi\equiv -\rho$, $scal_{\hat{g}}=0$ and the potential function $\tilde{h}$ can be expressed in the form $\tilde{h}(x)=-\frac{\rho}{2}|x|^{2}+\langle b,x\rangle+c$ for some $b\in\mathbb{R}^{n+d}$ and $c\in\mathbb{R}$; while if $\psi(x_{0})=\psi_{*}=0$ for some $x_{0}\in B$, then $\psi\equiv 0$, the scalar curvature $scal_{\hat{g}}$ is positive constant and $\nabla \tilde{h}$ is a Killing vector field.

		\item If $\rho=0$, then $R_{*}=0$. Furthermore, if $\psi(x_{0})=\psi_{*}=0$ for some $x_{0}\in B$, then $\psi\equiv0$, $scal_{\hat{g}}=0$ and $\nabla \tilde{h}$ is a Killing vector field.

		\item If $\rho<0$, then $-\rho\leq R_{*}\leq 0$.  Furthermore, if $\psi(x_{0})=\psi_{*}=0$ for some $x_{0}\in B$, then $\psi\equiv 0$, the scalar curvature $scal_{\hat{g}}$ is negative constant and $\nabla \tilde{h}$ is a Killing vector field; while if $\psi(x_{0})=\psi_{*}=-\rho$ for some $x_{0}\in B$, then $\psi\equiv -\rho$, $scal_{\hat{g}}=0$ and the potential function $\tilde{h}$ can be expressed in the form $\tilde{h}(x)=-\frac{\rho}{2}|x|^{2}+\langle b,x\rangle+c$ for some $b\in\mathbb{R}^{n+d}$ and $c\in\mathbb{R}$.

	\end{enumerate}
	
\end{theorem}

In addition we obtain

\begin{theorem}\label{th2}
	Under the same hypotheses of Theorem \ref{th1}, we have
	\begin{enumerate}
		\item If $\rho>0$, then $\tilde{h}(x)\geq-\frac{\rho}{2}r(x)^{2}+c_{1}r(x)+c_{2}$,\vspace{0,1cm}
		\item If $\rho\leq0$, then $\tilde{h}(x)\geq c_{1}r(x)+c_{2}$,
	\end{enumerate}
	for some fixed constants $c_{1}$, $c_{2}$, where $r(x)$ is the distance function from some point $p\in\hat{M}$.
\end{theorem}

In order to proceed we recall an important result which gives a characterization for warped product gradient Yamabe soliton.

\begin{proposition}\label{lemma}\textup{(\cite{tokura2019warped})}
	$(B^{n}\times_{f}F^{d},\hat{g},\nabla \tilde{h},\rho)$ is a \textbf{gradient Yamabe soliton} if, and only if, $(B^{n},g_{B},\nabla h,\lambda)$ is an \textbf{almost gradient Yamabe soliton} with soliton function 
	\begin{equation}\label{l1}
	\lambda=-\frac{\lambda_{F}}{f^2}+\frac{2d}{f}\Delta f+\frac{d(d-1)}{f^{2}}g_{B}(\nabla f,\nabla f)+\rho,
	\end{equation}
	and scalar curvature
	\begin{equation}\label{l2}
	scal_{g_B}=\frac{1}{f}g_{B}(\nabla f,\nabla h)+\lambda.
	\end{equation}
\end{proposition}

The previous result shows that the estimates for warping function might be applied in the study of gradient Yamabe solitons warped products. For instance, by the strong maximum principle, all expanding gradient Yamabe soliton $\mathbb{H}^{n}\times_{f}\mathbb{R}^{d}$ with $\rho=-n(n-1)$ are standard Riemannian product if $f$ attains  its maximum.

Consider the change $f=v^{\frac{2}{d+1}}$, then equation \eqref{l1} and \eqref{l2} turns out to be
\begin{equation}\label{9'}
\Delta_{w}v-\frac{1}{dp}(scal_{g_B}-\rho)v-\frac{1}{dp}\lambda_{F}v^{1-p}=0,
\end{equation}
 where \[w=-\frac{h}{2d},\quad p=\frac{4}{d+1},\]
 and $\Delta_{w}=e^{w}div(e^{-w}\nabla w)$, is the so called drifting Laplacian on the Bakry-\'Emery geometry.
 
 The attainability of the maximum by the function is something intimately related to the behaviour of its gradient. An interesting question is that whether or not we have local gradient estimates for positive warping solutions to the above equation.

In order to do this, in the remainder of this paper, we focus our attention on gradient estimates for the positive solutions to the nonlinear equation \eqref{9'}. We mainly follow the means of P. Li and S.T. Yau's proof in \cite{li1986parabolic}.

\begin{theorem}\label{gradient}Let $(B^{n}\times_{f}F^{d},\hat{g},\nabla \tilde{h},\rho)$ be a complete gradient Yamabe soliton satisfying
	$$^{B}Ric+Hessw-\frac{1}{m}dw\otimes dw\geq -K,\quad\Delta_{w}scal_{g_B}\leq\theta,\quad|\nabla scal_{g_B}|\leq\gamma,\quad\text{where}\quad w=-\frac{h}{2d},$$
	for $K\geq0$ in the metric ball $B(p,2R)$ of the base. Then, for any $\beta\in(0,1)$, the warping function $f$ satisfies the following gradient estimates:
	
	\begin{itemize}
		\item If $\lambda_{F}<0$, we have
		
		$$\beta\frac{|\nabla f|^{2}}{f^{2}}-\frac{scal_{B}-\rho}{d(d+1)}-\frac{\lambda_{F}}{d(d+1)f^{2}}\leq \frac{4(n+m)}{\beta(d+1)^{2}}B+4\sqrt{\frac{(n+m)C}{2\beta(d+1)^{4}}},\hspace{2,8cm}$$
		
		\item If $\lambda_{F}\geq0$, assume that $f$ is bounded in $B(p,2R)$, then we have
		
		$$\beta\frac{|\nabla f|^{2}}{f^{2}}-\frac{scal_{B}-\rho}{d(d+1)}-\frac{\lambda_{F}}{d(d+1)f^{2}}\leq \frac{4(n+m)}{\beta(d+1)^{2}}\left[B+\frac{\lambda_{F}M}{d}\right]+4\sqrt{\frac{(n+m)(C+D)}{2\beta(d+1)^{4}}},$$
		
	\end{itemize}
	in $B(p,R)$, where 
	
	\begin{eqnarray*}
		M&=&\displaystyle\sup_{B(p,2R)}f^{-2},\nonumber\\
		B&=&\frac{(n+m)c_{1}^{2}}{4R^{2}\beta(1-\beta)}+\frac{(n-1+R\sqrt{nK})c_{1}+c_{2}+2c_{1}^{2}}{R^2},\nonumber\\
		C&=&\frac{3\beta}{2} \left[\frac{n+m}{4}\left(\frac{\gamma}{dp}\right)^{4}\frac{(1-\beta)^{2}}{\beta^{4}}\varepsilon^{-1}\right]^{\frac{1}{3}}+\frac{\beta(n+m)}{2}(1-\varepsilon)^{-1}(1-\beta)^{-2}K^{2}+\frac{\theta(d+1)}{4d}, \nonumber\\
		D&=&\frac{\beta(n+m)}{2(1-\varepsilon)(1-\beta)^{2}}\Bigg{[}\left(\frac{(d+1)(\beta+\frac{4}{d+1}+1)\lambda_{F}M}{8d}\right)^{2}
		+\frac{(d+1)(\beta+\frac{4}{d+1}+1)K\lambda_{F}M}{4d}\Bigg{]},\nonumber
	\end{eqnarray*}
	$c_{1}$, $c_{2}$ are positive constants and $\varepsilon\in(0,1)$.
\end{theorem}

Letting $R\to\infty$ we get the following global gradient estimates.

\begin{cor}Let $(B^{n}\times_{f}F^{d},\hat{g},\nabla \tilde{h},\rho)$ be a complete gradient Yamabe soliton satisfying
	$$^{B}Ric+Hessw-\frac{1}{m}dw\otimes dw\geq -K,\quad\Delta_{w}scal_{g_B}\leq\theta,\quad|\nabla scal_{g_B}|\leq\gamma,\quad\text{where}\quad w=-\frac{h}{2d}.$$
	for $K\geq0$ in the base $B^{n}$. Then for any $\beta\in(0,1)$, the warping function $f$ satisfies the following gradient estimate
	
	\begin{itemize}
		\item If $\lambda_{F}<0$, we have
		
		$$\beta\frac{|\nabla f|^{2}}{f^{2}}-\frac{scal_{B}-\rho}{d(d+1)}-\frac{\lambda_{F}}{d(d+1)f^{2}}\leq 4\sqrt{\frac{(n+m)C}{2\beta(d+1)^{4}}},\hspace{3,7cm}$$
		
		\item If $\lambda_{F}=0$, we have
		
		$$\beta\frac{|\nabla f|^{2}}{f^{2}}-\frac{scal_{B}-\rho}{d(d+1)}\leq 4\sqrt{\frac{(n+m)C}{2\beta(d+1)^{4}}},\hspace{5.7cm}$$
		
		\item If $\lambda_{F}>0$, assume that $f$ is bounded, then we have
		
		$$\beta\frac{|\nabla f|^{2}}{f^{2}}-\frac{scal_{B}-\rho}{d(d+1)}-\frac{\lambda_{F}}{d(d+1)f^{2}}\leq \frac{4M'\lambda_{F}(n+m)}{d\beta(d+1)^{2}}+4\sqrt{\frac{(n+m)(C+D)}{2\beta(d+1)^{4}}},$$
		
	\end{itemize}
	in $B^{n}$, where 
	
	\begin{eqnarray*}
		M'&=&\displaystyle\sup_{B}f^{-2},\nonumber\\
		C&=&\frac{3\beta}{2} \left[\frac{n+m}{4}\left(\frac{\gamma}{dp}\right)^{4}\frac{(1-\beta)^{2}}{\beta^{4}}\varepsilon^{-1}\right]^{\frac{1}{3}}+\frac{\beta(n+m)}{2}(1-\varepsilon)^{-1}(1-\beta)^{-2}K^{2}+\frac{\theta(d+1)}{4d}, \nonumber\\
		D&=&\frac{\beta(n+m)}{2(1-\varepsilon)(1-\beta)^{2}}\Bigg{[}\left(\frac{(d+1)(\beta+\frac{4}{d+1}+1)\lambda_{F}M}{8d}\right)^{2}
		+\frac{(d+1)(\beta+\frac{4}{d+1}+1)K\lambda_{F}M}{4d}\Bigg{]},\nonumber
	\end{eqnarray*}
and $\varepsilon\in(0,1)$.
\end{cor}

As an application, we obtain the following two results.

\begin{cor}\label{cor}There is no complete gradient Yamabe soliton $(B^{n}\times_{f}F^{d},\hat{g},\nabla \tilde{h},\rho)$ with 
$$^{B}Ric+Hessw-\frac{1}{m}dw\otimes dw\geq 0,\quad scal_{B}=\text{constant}\leq\rho,\quad \lambda_{F}<0,\quad \text{where}\quad w=-\frac{h}{2d}.$$
\end{cor}

\begin{cor}\label{cor2}There is no complete gradient Yamabe soliton $(B^{n}\times_{f}F^{d},\hat{g},\nabla \tilde{h},\rho)$ with 
	$$^{B}Ric+Hessw-\frac{1}{m}dw\otimes dw\geq 0,\quad scal_{B}=\text{constant}<\rho,\quad \lambda_{F}=0,\quad \text{where}\quad w=-\frac{h}{2d}.$$
\end{cor}

\begin{observation} Corollary
	 \ref{cor} and Corollary \ref{cor2} produce constraints on the construction of Yamabe solitons. For instance, consider $\hat{M}=\mathbb{R}\times_{f}\mathbb{H}^{d}$, then there does not exist complete shrinking or steady gradient Yamabe soliton metric on $\hat{M}$, with potential function satisfying
	
	$$h''+\frac{h'^{2}}{2d}\leq0$$

In the trivial case, the above example bring to light that the product manifold $\mathbb{R}^{n}\times_{f}\mathbb{H}^{d}$ does not admit a complete metric with constant positive scalar curvature.
\end{observation}

\section{Examples}
\label{examples}
In this section we present some examples of warped product gradient Yamabe solitons. These examples are interesting to guide our intuition about general properties of the gradient yamabe solitons. First denote
$$\mathbb{R}_{+}^{3}=\{(x_{1},x_{2},x_{3})\in\mathbb{R}^{3};x_{3}>0\},\quad \mathbb{R}_{\ast}^{3}=\{(x_{1},x_{2},x_{3})\in\mathbb{R}^{3};x_{1}+x_{2}+x_{3}>0\},$$
and 
$$ds^2=(dx_{1})^{2}+(dx_{2})^{2}+(dx_{3})^{2}$$

\begin{example}(Einstein Warped product) Consider the product manifold $\hat{M}=\mathbb{R}_{+}^{3}\times \mathbb{R}^{3}$ furnished with metric $\hat{g}=ds_{1}^2+f^2 ds$, where,
	$$ds_{1}^2=\frac{ds^2}{x_{3}^{2}},\quad f(x_{1},x_{2},x_{3})=\frac{1}{x_{3}}.$$
	A straightforward computation shows that $\mathbb{H}^{3}\times _{f}\mathbb{R}^{3}$ is a complete noncompact Einstein warped product with $^{\hat{M}}Ric=-\frac{5}{6}\hat{g}$. Then $\hat{M}$ is trivial gradient Yamabe soliton.
\end{example}

\begin{example}Consider the product manifold $\hat{M}=\mathbb{R}_{\ast}^{3}\times \mathbb{R}^3$ furnished with metric $\hat{g}=ds_{1}^2+f^2 ds$, where,
	$$ds_{1}^2=20\frac{ds^2}{x_{1}+x_{2}+x_{3}},\quad f(x_{1},x_{2},x_{3})=\sqrt{\frac{20}{x_{1}+x_{2}+x_{3}}}.$$
A straightforward computation shows that $\hat{M}$ is a steady gradient Yamabe soliton with potential function $$h(x_{1},x_{2},x_{3})=20\log(x_{1}+x_{2}+x_{3}).$$
\end{example}
The next example concerns a complete noncompact gradient Yamabe soliton.
\begin{example}\label{example1}Let $(\mathbb{S}^3,dr^2)$ the standard 3-dimensional sphere and consider the product manifold $\hat{M}=\mathbb{R}_{+}^{3}\times(\mathbb{S}^3\times\mathbb{R}_{+}^3)$ furnished with metric $\hat{g}=ds_{1}^2+f^2(dr^2+ds_{1}^2)$, where,
	$$ds_{1}^2=\frac{ds^2}{x_{3}^{2}},\quad f(x_{1},x_{2},x_{3})=\sqrt[7]{x_{3}^{4}}.$$	
	A straightforward computation shows that $\mathbb{H}^{3}\times _{f}(\mathbb{S}^{3}\times \mathbb{H}^3)$ is a complete noncompact trivial gradient Yamabe soliton with $\rho=-6$.
\end{example}

\section{Proofs}
\label{proofs}

\begin{myproof}{Theorem}{\ref{th1}} 

We know that a gradient Yamabe soliton $(M^{n},g,\nabla h, \rho)$ with scalar curvature $scal_{g}$, satisfy (cf. \cite{barbosa2013conformal}):

\begin{equation}\label{1}(n-1)\Delta scal_{g}+\frac{1}{2}g(\nabla scal_{g},\nabla h)+scal_{g}(scal_{g}-\rho)=0.
\end{equation}

From Proposition \ref{lemma}, we obtain:
\begin{equation}\label{10}
\begin{split}
scal_{\hat{g}}&=\left[scal_{g_B}+\frac{\lambda_{F}}{f^{2}}-2d\frac{\Delta f}{f}-\frac{d(d-1)}{f^{2}}g_{B}(\nabla f,\nabla f)\right]\circ\pi,\\
&=[scal_{g_B}+\rho-\lambda]\circ\pi.
\end{split}
\end{equation}

Since $(\hat{M}^{n+d},\hat{g},\nabla\tilde{h},\rho)$ is a gradient Yamabe soliton, combining \eqref{1} and \eqref{10} we get:

\begin{equation*}
\begin{split}
(n+d-1)\Delta (scal_{g_B}-\lambda)+\frac{d(n+d-1)}{f}g_{B}(\nabla(scal_{g_{B}}-\lambda),\nabla f)+\\+\frac{1}{2}g_{B}(\nabla (scal_{g_B}-\lambda),\nabla h)
+(scal_{g_{B}}-\lambda+\rho)(scal_{g_{B}}-\lambda)=0.
\end{split}
\end{equation*}
Therefore,
\begin{equation}\label{3'}
\Delta_{w}\psi=-\frac{1}{n+d-1}(\psi+\rho)\psi,
\end{equation}
where $\psi=scal_{g_B}-\lambda$ and $w=-d \log f-2^{-1}(n+d-1)^{-1}h$.

We proceed observing that, by Proposition 3.3 of \cite{pigola2010ricci}, the following volume estimate holds
\begin{equation}\label{5}
Vol_{w}(B_{r})\leq Ae^{Br^{2}},\qquad A,B\in(0,\infty).
\end{equation}
In particular,
\begin{equation*}
\liminf_{r\to\infty}\frac{\log Vol_{w}(B_{r})}{r^2}\leq C<\infty,\qquad C\in\mathbb{R}.
\end{equation*}
From equation \eqref{3'}, setting $\psi_{-}=\max\{-\psi,0\}$ we immediately deduce that
\begin{equation*}
\Delta_{w}\psi_{-}=\frac{\psi_{-}^{2}-\rho\psi_{-}}{n+d-1}.
\end{equation*} 
Then, applying Theorem 12 of \cite{pigola2011remarks} with the choices
$$a(x)=-\rho(n+d-1)^{-1},\quad b(x)=(n+d-1)^{-1},\quad\sigma =2,$$
we obtain that $\psi_{-}$ is bounded from above, or equivalently, 

\begin{equation*}\psi_{\ast}=\inf_{B}\psi>-\infty.
\end{equation*}

Next, again by \eqref{5}, the weak minimum principle at infinity for $\Delta_{w}$ holds (see Theorem 9 of \cite{pigola2011remarks}). Then there exist a sequence $\{x_{k}\}$ such that,
 $$\psi(x_{k})\to\psi_{*},\qquad\Delta_{w}\psi(x_{k})\geq-\frac{1}{k},$$
and taking the limit in \eqref{3'} along $\{x_{k}\}$ we obtain
\begin{equation}\label{quadratic}
 (\psi_{*}+\rho)\psi_{*}\leq0.
\end{equation}
However $\psi\circ\pi=scal_{\hat{g}}-\rho$, so that the claimed bounds on $R_{\ast}=\inf_{\hat{M}}scal_{\hat{g}}$ in the statement of Theorem \ref{th1} follow immediately from \eqref{quadratic}.

\textbf{Case I: $\rho>0$}. Assume that $\psi(x_{0})=\psi_{*}=-\rho$, for some $x_{0}\in B$. Then by \eqref{3'} the  non negative function $l(x)=\psi(x)+\rho$ satisfies
\begin{equation}\label{mp1}
\Delta_{w}l-\frac{\rho}{n+d-1}l=-\frac{l^2}{n+d-1}\leq0.
\end{equation}
We let
\begin{equation*}
\Omega_{0}:=\{x\in B\hspace{0.2cm} ;\hspace{0.2cm} l(x)=0\},
\end{equation*}
$\Omega_{0}$ is closed and nonempty since $x_{0}\in \Omega_{0}$. Let now $y\in \Omega_{0}$, then applying the maximum principle (see \cite{gilbarg2015elliptic} p. 35 ) to \eqref{mp1}, we obtain that $l(x)=0$ in a neighborhood of $y$, so that $\Omega_{0}$ is open. Connectedness of $B$ yields $\Omega_{0}=B$. Therefore $\psi\equiv-\rho$ or, equivalently,

$$scal_{\hat{g}}=0.$$
Combining $scal_{\hat{g}}=0$ with equation \eqref{def2} gives
\begin{equation}\label{hess}
-\rho \hat{g}=Hess_{\hat{g}} \tilde{h}.
\end{equation}
Thus, by Theorem 1 of \cite{pigola2011remarks}, $\hat{M}$ is isometric to $\mathbb{R}^{n+d}$ and, solving equation \eqref{hess}, we get

$$\tilde{h}(x)=-\frac{\rho}{2}|x|^{2}+\langle b,x\rangle+c,$$
for some $b\in\mathbb{R}^{n+d}$ and $c\in\mathbb{R}$,
which proves the first assertion of item $(1)$.

Analogously, if $\psi(x_{0})=\psi_{*}=0$, for some $x_{0}\in B$, we deduce that $\psi\circ\pi=scal_{\hat{g}}-\rho\equiv0$ and therefore, $scal_{\hat{g}}$ is positive constant and  $\nabla \tilde{h}$ is a Killing vector field.

\textbf{Case II: $\rho=0$}. Assume that $\psi(x_{0})=\psi_{*}=0$ for some $x_{0}\in B$. Then by \eqref{3'} the non negative function $\psi(x)$ satisfies
\begin{equation*}
\Delta_{w}\psi=-\frac{1}{n+d-1}\psi^2\leq0.
\end{equation*}
By the maximum principle we conclude that $\psi\circ\pi=scal_{\hat{g}}\equiv0$, and therefore $\nabla\tilde{h}$ is a Killing vector field, thus concluding the proof of item $(2)$.

\textbf{Case III: $\rho<0$}. Assume that $\psi(x_{0})=\psi_{*}=0$ for some $x_{0}\in B$. From \eqref{3'}, we have
\begin{equation*}
\Delta_{w}\psi+\frac{\rho}{n+d-1}\psi=-\frac{\psi^2}{n+d-1}\leq0.
\end{equation*}
Since $\psi(x)\geq\psi_{\ast}=0$, by the maximum principle we conclude $\psi\circ\pi=scal_{\hat{g}}-\rho\equiv0$ and therefore, $scal_{\hat{g}}$ is negative constant and $\nabla\tilde{h}$ is a Killing vector field, which proves the first assertion of item $(3)$.

Analogously, suppose that $\psi(x_{0})=\psi_{*}=-\rho$ for some $x_{0}\in B$. Then,  again by the maximum principle $\psi\equiv-\rho$ or, what we already know to be the same as the following
$$scal_{\hat{g}}=0.$$
Combining $scal_{\hat{g}}=0$ with equation \eqref{def2} we obtain
\begin{equation}\label{hess3}
-\rho \hat{g}=Hess_{\hat{g}}\tilde{h}.
\end{equation}
Thus, by Theorem 1 of \cite{pigola2011remarks}, $\hat{M}$ is isometric to $\mathbb{R}^{n+d}$ and, solving equation \eqref{hess3}, we get

$$\tilde{h}(x)=-\frac{\rho}{2}|x|^{2}+\langle b,x\rangle+c,$$
for some $b\in\mathbb{R}^{n+d}$ and $c\in\mathbb{R}$.

\end{myproof}

\begin{myproof}{Theorem}{\ref{th2}}Let $r(x):=r(x,x_{0})$ the distance function from a fixed point $x_{0}$ and consider $\alpha(s):[0,r]\rightarrow \hat{M}$ a minimizing geodesic emanated from $x_{0}=\alpha(0)$. Then
\begin{eqnarray*}
\frac{d}{ds}\Big{|}_{r}\tilde{h}(\alpha(s))&=&\hat{g}(\nabla \tilde{h}, \alpha'(r)),\\
&=&\int_{0}^{r}\frac{d}{ds}\hat{g}(\nabla\tilde{h},\alpha'(s))ds+\hat{g}(\nabla \tilde{h}, \alpha'(0)),\\
&=&\int_{0}^{r}Hess_{\hat{g}} \tilde{h}(\alpha',\alpha')ds+\hat{g}(\nabla\tilde{h},\alpha'(0)).
\end{eqnarray*}

By Theorem \ref{th1} we have that

\begin{equation*}
Hess_{\hat{g}}\tilde{h}(\alpha',\alpha') =scal_{\hat{g}}-\rho\geq \begin{cases} -\rho &\text{if } \rho>0 \\[10pt]
0 & \text{if } \rho\leq0 
\end{cases}
\end{equation*}
Therefore,
\begin{equation*}
\frac{d}{ds}\Big{|}_{r}\widetilde{h}(\alpha(s))\geq\begin{cases} -\rho r+\hat{g}(\nabla\tilde{h},\alpha'(0)) &\text{if } \rho>0 \\[10pt]
\hat{g}(\nabla\tilde{h},\alpha'(0)) & \text{if } \rho\leq0 
\end{cases}
\end{equation*}
Integrating the above inequalities along $\alpha(s)$ yields Theorem \ref{th2}.

\end{myproof}

\begin{myproof}{Theorem}{\ref{gradient}}
Consider the change $f=v^{\frac{2}{d+1}}$. Then by Proposition \ref{lemma} we have that
\begin{equation}\label{9}
\Delta_{w}v-\frac{1}{dp}(scal_{g_B}-\rho)v-\frac{1}{dp}\lambda_{F}v^{1-p}=0,\quad \text{where}\quad w=-\frac{h}{2d},\quad p=\frac{4}{d+1}.
\end{equation}
Write $g_{B}=\langle\cdot,\cdot\rangle=|\cdot|^{2}$ for simplicity. Let $v$ a positive solution to \eqref{9}, then $u=\log v$ satisfies the equation

\begin{equation*}
\Delta_{w}u=(\beta-1)|\nabla u|^{2}-L,\quad \text{where}\quad L=\beta|\nabla u|^{2}-\frac{1}{dp}(scal_{g_B}-\rho)-\frac{\lambda_{F}}{dp}e^{-pu},\quad \beta\in(0,1).
\end{equation*}

Now, consider a cut-off function $\xi$ satisfying

\begin{equation*}
\xi(r)= \begin{cases} 1 &\text{if } r\in[0,1] \\[10pt]
0 & \text{if } r\in[2,\infty) 
\end{cases},\qquad -c_{1}\leq\frac{\xi'(r)}{\xi^{\frac{1}{2}}(r)}\leq0,\qquad -c_{2}\leq\xi''(r),\qquad c_{1},c_{2}\in (0,\infty),
\end{equation*}
and define
\begin{equation*}
\psi(x)=\xi\left(\frac{r(x)}{R}\right),
\end{equation*}
where $r(x)$ is the distance function starting from $p$ to $x$.
Using an argument of Calabi \cite{calabi1958extension}( see also Cheng and Yau \cite{cheng1975differential}), we can assume without loss of generality that the function $\psi$ is smooth in $B(p,2R)$. Then, the function defined by $G=\psi L$ is smooth in $B(p,2R)$. 

Let $x_{0}\in B(p,2R)$ be a point at which $G$ attains its maximum value $G_{max}$, and suppose that $G_{max}>0$ (otherwise the proof is trivial). At the point $x_{0}$, we have
\begin{equation*}
\nabla (G)=\psi\nabla L+L\nabla \psi=0.
\end{equation*}
Moreover,
\begin{equation}\label{general}
\begin{split}
0&\geq\Delta_{w}G,\\
&=\psi\Delta_{w}L+L\Delta_{w}\psi+2\langle\nabla \psi,\nabla L\rangle,\\
&=\psi\Delta_{w}L+L\Delta_{w}\psi-2L\frac{|\nabla\psi|^{2}}{\psi}.
\end{split}
\end{equation}
In order to estimate the left side of \eqref{general} we prove the following lemma:

\begin{lema}\label{lemma2}Let $(B^{n},g_{B})$ be a complete noncompact Riemannian  manifold satisfying
	
\begin{equation}\label{ricci hessian}^{B}Ric+Hessw-\frac{1}{m}dw\otimes dw\geq -K,\quad\text{where}\quad w=-\frac{h}{2d},
\end{equation}
for $K\geq0$ in the metric ball $B(p,2R)\subset B$, and let $L$ and $\psi$ as above. Then, we have

	\begin{equation}\label{a2'}
	\frac{|\nabla\psi|^{2}}{\psi}\leq\frac{c_{1}^{2}}{R^{2}},
	\end{equation}
	
	\begin{equation}\label{a3'}
	\Delta_{w}\psi\geq-\frac{(n-1+R\sqrt{nK})c_{1}+c_{2}}{R^2},
	\end{equation}
	
		\begin{equation}\label{a1}
		\begin{split}
	\Delta_{w}L\geq 2\beta\frac{(\Delta_{w}u)^{2}}{n+m}+\frac{2(\beta-1)}{dp}\langle\nabla u,\nabla scal_{g_B}\rangle-\frac{2\beta\lambda_{F}}{d}e^{-pu}|\nabla u|^{2}-2\langle\nabla u,\nabla L\rangle+\\-2\beta K|\nabla u|^{2}
	-\frac{\Delta_{w}scal_{g_B}}{dp}-\frac{\lambda_{F}}{d}e^{-pu}\left[(p-\beta+1)|\nabla u|^{2}+L\right].
	\end{split}
	\end{equation}
\end{lema}

\begin{myproof}{Lemma}{\ref{lemma2}}
Equation \eqref{a2'} follows from the calculation

\begin{equation*}
\frac{|\nabla\psi|^{2}}{\psi}=\frac{1}{\xi}\left\langle \xi'\frac{\nabla r}{R}, \xi'\frac{\nabla r}{R}\right\rangle=\frac{(\xi')^2}{\xi}\frac{1}{R^2}\langle\nabla r,\nabla r\rangle\leq\frac{c_{1}^{2}}{R^{2}}.\\[10pt]
\end{equation*}

It has been shown by Qian \cite{qian1998comparison}, the following estimate

\begin{equation*}
\Delta_{w}r^{2}\leq n\left(1+\sqrt{1+\frac{4Kr^2}{n}}\right),
\end{equation*}
which implies
\begin{eqnarray}
\Delta_{w}r&=&\frac{1}{2r}\left(\Delta_{w}r^2-2|\nabla r|^{2}\right),\nonumber\\
&\leq&\frac{n-2}{2r}+\frac{n}{2r}\left(1+\sqrt{1+\frac{4Kr^2}{n}}\right),\nonumber\\
&=&\frac{n-1}{r}+\sqrt{nK}\nonumber.
\end{eqnarray}
Then, we obtain
\begin{equation*}
\Delta_{w}\psi=\frac{\xi''(r)|\nabla r|^{2}}{R^2}+\frac{\xi'(r)\Delta_{w}r}{R}\geq-\frac{(n-1+R\sqrt{nK})c_{1}+c_{2}}{R^2},
\end{equation*}
which proves \eqref{a3'}.

From the Bochner formula for the $m$-Bakry-\'Emery Ricci tensor and
the lower bound hypothesis \eqref{ricci hessian} we obtain

\begin{eqnarray*}
	\frac{1}{2}\Delta_{w}|\nabla u|^{2}\geq\frac{(\Delta_{w}u)^{2}}{n+m}+\langle\nabla u,\nabla\Delta_{w}u\rangle-K|\nabla u|^{2}.
\end{eqnarray*}
Therefore, 
\begin{equation*}
\begin{split}
\Delta_{w}L&=\beta\Delta_{w}|\nabla u|^{2}-\frac{1}{dp}\Delta_{w}scal_{g_B}-\frac{\lambda_{F}}{dp}\Delta_{w}e^{-pu}\\
&\geq 2\beta\frac{(\Delta_{w}u)^{2}}{n+m}+2\beta\langle\nabla u,\nabla\Delta_{w}u\rangle-2\beta K|\nabla u|^{2}-\frac{1}{dp}\Delta_{w}scal_{g_B}+\\
&-\frac{\lambda_{F}}{d}e^{-pu}(p|\nabla u|^{2}-\Delta_{w}u).
\end{split}
\end{equation*}
Notice that

\begin{eqnarray}
2\beta\langle\nabla u,\nabla\Delta_{w}u\rangle&=&2\beta\langle\nabla u,\nabla\left[\left(1-\frac{1}{\beta}\right)\left(\frac{1}{dp}(scal_{g_B}-\rho)+\frac{\lambda_{F}e^{-pu}}{dp}\right)-\frac{L}{\beta}\right]\rangle,\nonumber\\
&=&2\beta\left(1-\frac{1}{\beta}\right)\frac{1}{dp}\langle\nabla u,\nabla scal_{g_B}\rangle+2\beta\frac{\lambda_{F}}{dp}\langle\nabla u,\nabla e^{-pu}\rangle-2\langle\nabla u,\nabla L\rangle,\nonumber\\
&=&2(\beta-1)\frac{1}{dp}\langle\nabla u,\nabla scal_{g_B}\rangle-\frac{2\beta\lambda_{F}}{d}e^{-pu}|\nabla u|^{2}-2\langle\nabla u,\nabla L\rangle,\nonumber
\end{eqnarray}
and

\begin{eqnarray}
\frac{\lambda_{F}}{d}e^{-pu}(p|\nabla u|^{2}-\Delta_{w}u)&=&\frac{\lambda_{F}}{d}e^{-pu}\left[p|\nabla u|^{2}-(\beta-1)|\nabla u|^{2}+L\right],\nonumber\\
&=&\frac{\lambda_{F}}{d}e^{-pu}\left[(p-\beta+1)|\nabla u|^{2}+L\right].\nonumber
\end{eqnarray}
It follows that

\begin{equation*}
\begin{split}
\Delta_{w}L\geq 2\beta\frac{(\Delta_{w}u)^{2}}{n+m}+\frac{2(\beta-1)}{dp}\langle\nabla u,\nabla scal_{g_B}\rangle-\frac{2\beta\lambda_{F}}{d}e^{-pu}|\nabla u|^{2}-2\langle\nabla u,\nabla L\rangle+\\-2\beta K|\nabla u|^{2}
-\frac{\Delta_{w}scal_{g_B}}{dp}-\frac{\lambda_{F}}{d}e^{-pu}\left[(p-\beta+1)|\nabla u|^{2}+L\right].
\end{split}
\end{equation*}
which completes the proof of lemma.
\end{myproof}\\[1pt]

Proceeding, using Lemma \ref{lemma2} and \eqref{general}, we obtain at the point $x_{0}$,

\begin{eqnarray*}
\psi\Bigg{(}2\beta\frac{(\Delta_{w}u)^{2}}{n+m}+2(\beta-1)\frac{1}{dp}\langle\nabla u,\nabla scal_{g_B}\rangle-\frac{2\beta\lambda_{F}}{d}e^{-pu}|\nabla u|^{2}-2\langle\nabla u,\nabla L\rangle-2\beta K|\nabla u|^{2}+\\
-\frac{\Delta_{w}scal_{g_B}}{dp}-\frac{\lambda_{F}}{d}e^{-pu}\left[(p-\beta+1)|\nabla u|^{2}+L\right]\Bigg{)}\leq LH,
\end{eqnarray*}
where

\begin{equation*}
H=\left(\frac{(n-1+R\sqrt{nK})c_{1}+c_{2}+2c_{1}^{2}}{R^2}\right).
\end{equation*}

From the fact that $0\leq\psi\leq1$, we have
$$
-2\psi\langle\nabla u,\nabla L\rangle=2L\langle\nabla u, \nabla \psi\rangle\geq-2L|\nabla u||\nabla\psi|\geq-\frac{2c_{1}}{R}\psi^{\frac{1}{2}}L|\nabla u|.
$$
Then

\begin{equation}\label{final}
\begin{split}
	2\beta\psi\frac{(\Delta_{w}u)^{2}}{n+m}+2(\beta-1)\frac{\psi}{dp}\langle\nabla u,\nabla scal_{g_B}\rangle-\frac{(\beta+p+1)\lambda_{F}\psi}{d}e^{-pu}|\nabla u|^{2}+\\
	-\frac{2c_{1}}{R}\psi^{\frac{1}{2}}L|\nabla u|
	-2\beta\psi K|\nabla u|^{2}
	-\psi\frac{\Delta_{w}scal_{g_B}}{dp}-\frac{\lambda_{F}\psi}{d}e^{-pu}L\leq LH.
\end{split}
\end{equation}
In the sequel, we distinguish between two cases: (a) $\lambda_{F}<0$ and (b) $\lambda_{F}\geq0$.

Case (a): $\lambda_{F}<0$. Since

$$\Delta_{w}scal_{g_B}\leq\theta(2R),\qquad|\nabla scal_{g_B}|\leq\gamma(2R),$$
then \eqref{final} yields

\begin{eqnarray*}
	2\beta\psi\frac{(\Delta_{w}u)^{2}}{n+m}+2(\beta-1)\frac{\psi}{dp}\gamma|\nabla u|-\frac{2c_{1}}{R}\psi^{\frac{1}{2}}L|\nabla u|
	-2\beta\psi K|\nabla u|^{2}-\frac{\theta\psi}{dp}\leq LH.
\end{eqnarray*}
Multiplying both sides of the above equation by $\psi$ and using the fact that $0\leq\psi\leq 1$, we obtain

\begin{eqnarray*}
	2\beta\frac{(\psi\Delta_{w}u)^{2}}{n+m}+2(\beta-1)\frac{\psi^{\frac{1}{2}}}{dp}\gamma|\nabla u|-\frac{2c_{1}}{R}\psi^{\frac{3}{2}}L|\nabla u|
	-2\beta\psi K|\nabla u|^{2}-\frac{\theta}{dp}\leq \psi LH.
\end{eqnarray*}
Let

$$y=\psi|\nabla u|^{2},\qquad z=\psi\left(\frac{1}{dp}(scal_{g_B}-\rho)+\frac{\lambda_{F}}{dp}e^{-pu}\right).$$
Then we have

\begin{eqnarray*}
	\frac{2\beta}{n+m}\Big\{ (y-z)^{2}+\frac{(\beta-1)\gamma(n+m)y^{\frac{1}{2}}}{\beta dp}-\frac{(n+m)c_{1}}{R}y^{\frac{1}{2}}\left(y-\frac{z}{\beta}\right)-(n+m)Ky\Big\}\leq\\
	\leq \psi LH+\frac{\theta}{dp}.
\end{eqnarray*}
From Li-Yau's arguments (\cite{li1986parabolic}, pg.161-162), for any $0<\varepsilon<1$ we obtain that
\begin{eqnarray*}
	\frac{2\beta}{n+m}\Bigg\{(\beta y-z)^{2}-\frac{(n+m)^{2}c_{1}^{2}}{8R^{2}\beta^{2}(1-\beta)}(\beta y-z)-\frac{3}{4}4^{-\frac{1}{3}}(n+m)^{\frac{4}{3}}\left[\left(\frac{\gamma}{dp}\right)^{4}\frac{(1-\beta)^{2}}{\beta^{4}}\varepsilon^{-1}\right]^{\frac{1}{3}}+\nonumber\\
	-\frac{(n+m)^{2}}{4}(1-\varepsilon)^{-1}(1-\beta)^{-2}K^{2}\Bigg\}
	\leq (\beta y-z)H+\frac{\theta}{dp}.
\end{eqnarray*}
Hence,
\begin{equation*}
\frac{2\beta}{n+m}(\psi L)^{2}-B(\psi L)-C\leq0,
\end{equation*}
where,

\begin{eqnarray*}
B&=&\frac{(n+m)c_{1}^{2}}{4R^{2}\beta(1-\beta)}+H,\nonumber\\
C&=&\frac{3\beta}{2} \left[\frac{n+m}{4}\left(\frac{\gamma}{dp}\right)^{4}\frac{(1-\beta)^{2}}{\beta^{4}}\varepsilon^{-1}\right]^{\frac{1}{3}}+\frac{\beta(n+m)}{2}(1-\varepsilon)^{-1}(1-\beta)^{-2}K^{2}+\frac{\theta}{dp}. \nonumber
\end{eqnarray*}

Using the inequality $Az^{2}-Bz\leq C$, one obtain $z\leq\frac{2B}{A}+\sqrt{\frac{C}{A}}$. Then

\begin{eqnarray*}
\sup_{x\in B(p,R)}L(x)\leq(\psi L)(x_{0})\leq \frac{n+m}{\beta}B+\sqrt{\frac{n+m}{2\beta}}C^{\frac{1}{2}},
\end{eqnarray*}
and hence,

$$\beta|\nabla u|^{2}-\frac{1}{dp}(scal_{g_B}-\rho)-\frac{\lambda_{F}}{dp}e^{-pu}\leq \frac{n+m}{\beta}B+\sqrt{\frac{n+m}{2\beta}}C^{\frac{1}{2}}.$$

Replacing the function $u=\log f^{\frac{d+1}{2}}$ back into the above equation we obtain the desired estimate for $\lambda_{F}<0$.

Case (b): $\lambda_{F}\geq0$. Since

$$\Delta_{w}scal_{g_B}\leq\theta(2R),\qquad|\nabla scal_{g_B}|\leq\gamma(2R),$$
then \eqref{final} yields

\begin{eqnarray*}
	2\beta\psi\frac{(\Delta_{w}u)^{2}}{n+m}+2(\beta-1)\frac{\psi}{dp}\gamma|\nabla u|-\frac{(\beta+p+1)\lambda_{F}\psi}{d}e^{-pu}|\nabla u|^{2}-\frac{2c_{1}}{R}\psi^{\frac{1}{2}}L|\nabla u|+\\
	-2\beta\psi K|\nabla u|^{2}
	-\psi\frac{\theta}{dp}-\frac{\lambda_{F}\psi}{d}e^{-pu}L\leq LH.
\end{eqnarray*}
Multiplying both sides of the above equation by $\psi$, and using the fact that $0\leq\psi\leq 1$, we obtain

\begin{eqnarray*}
	2\beta\frac{(\psi\Delta_{w}u)^{2}}{n+m}+2(\beta-1)\frac{\psi^{\frac{1}{2}}}{dp}\gamma|\nabla u|-\frac{(\beta+p+1)\lambda_{F}\psi}{d}M|\nabla u|^{2}-\frac{2c_{1}}{R}\psi^{\frac{3}{2}}L|\nabla u|+\\
	-2\beta\psi K|\nabla u|^{2}
	-\psi\frac{\theta}{dp}-\frac{\lambda_{F}\psi}{d}ML\leq \psi LH,
\end{eqnarray*}
where $M=\sup_{B(p,2R)} e^{-pu}$.

Let

$$y=\psi|\nabla u|^{2},\qquad z=\psi\left(\frac{1}{dp}(scal_{g_B}-\rho)+\frac{\lambda_{F}}{dp}e^{-pu}\right).$$
Then we have

\begin{eqnarray*}
	\frac{2\beta}{n+m}\Bigg\{(y-z)^{2}-\frac{c_{1}(n+m)}{R}y^{\frac{1}{2}}\left(y-\frac{z}{\beta}\right)
	+(n+m)\frac{(\beta-1)}{\beta}\frac{y^{\frac{1}{2}}}{dp}\gamma+\nonumber\\-(n+m)\left[\frac{(\beta+p+1)\lambda_{F}}{2d\beta}M+K\right]y
	\Bigg\}\leq \psi L\left[H+\frac{\lambda_{F}M}{d}\right]+\frac{\theta}{dp}.
\end{eqnarray*}
Hence, we get that

\begin{eqnarray*}
	\frac{2\beta}{n+m}(\psi L)^{2}-\Bigg\{B+\frac{\lambda_{F}M}{d}\Bigg\}(\psi L)-\Bigg\{C+\frac{\beta(n+m)}{2(1-\varepsilon)(1-\beta)^{2}}\Bigg{[}\left(\frac{(\beta+p+1)\lambda_{F}M}{2dp}\right)^{2}+\nonumber\\
	+\frac{K(\beta+p+1)\lambda_{F}M}{dp}\Bigg{]}\Bigg\}\leq 0,
\end{eqnarray*}
and then

\begin{eqnarray*}
	\psi L\leq \frac{n+m}{\beta}\left[B+\frac{\lambda_{F}M}{d}\right]+\sqrt{\frac{n+m}{2\beta}}(C+D)^{\frac{1}{2}},
\end{eqnarray*}
where,
$$D=\frac{\beta(n+m)}{2(1-\varepsilon)(1-\beta)^{2}}\Bigg{[}\left(\frac{(\beta+p+1)\lambda_{F}M}{2dp}\right)^{2}
+\frac{K(\beta+p+1)\lambda_{F}M}{dp}\Bigg{]}.$$

Replacing the function $u=\log f^{\frac{d+1}{2}}$ back into the above equation we obtain the desired estimate for $\lambda_{F}\geq0$.
\end{myproof}

\end{document}